# On a generalization of $L^p$-differentiability

Daniel Spector

Department of Applied Mathematics, National Chiao Tung University, Hsinchu, Taiwan

October 14, 2015

## Abstract

In this paper we connect Calderón and Zygmund's notion of $L^p$-differentiability [5] with some recent characterizations of Sobolev spaces via the asymptotics of non-local functionals due to Bourgain, Brezis, and Mironescu [3]. We show how the results of the former can be generalized to the setting of the latter, while the latter results can be strengthened in the spirit of the former. As a consequence of these results we give several new characterizations of Sobolev spaces, a novel condition for whether a function of bounded variation is in the Sobolev space $W^{1,1}$, and complete the proof of a characterization of the Sobolev spaces recently claimed in [9, 10].

## 1 Introduction

### 1.1 $L^p$-differentiability and $L^p$-Taylor approximations

$L^p$-differentiability was introduced by Calderón and Zygmund in their study of the local properties of solutions of elliptic differential equations [5, 6]. It is a natural extension of classical differentiability in that it relaxes the requirement of the existence of a locally linear map in a uniform sense to its existence in an averaged sense. As the Sobolev spaces arise readily in the study of partial differential equations, it is not surprising that Sobolev functions possess an $L^p$-derivative. The following theorem asserting this fact was proven by Calderón and Zygmund [6, Theorem 12] (for a modern reference, see the monograph of Evans and Gariepy [8]).

**Theorem 1.1 (Calderón and Zygmund)** *Suppose* $1 \leq p < \infty$ *and* $f \in W^{1,p}(\mathbb{R}^N)$. *Then*

$$\left(\frac{1}{\epsilon^{pq}} \fint_{B(0,\epsilon)} |f(x+h) - f(x) - \nabla f(x) \cdot h|^{pq} \, dh\right)^{\frac{1}{pq}} \to 0 \qquad (1.1)$$

*for $\mathcal{L}^N$ almost every $x \in \mathbb{R}^N$, where $1 \leq q \leq \frac{N}{N-p}$ if $1 \leq p < N$, $1 \leq q < \infty$ if $p = N$, and $1 \leq q \leq \infty$ if $p > N$ (When $q = \infty$ the left side of (1.1) is understood to be $L^\infty_h(B(0,\epsilon))$ norm applied to the integrand).*

While $L^p$-differentiability is a necessary condition for the inclusion of a function in $W^{1,p}$, the fact that it does not characterize the space is readily seen by the improvement of exponent in Theorem 1.1. A natural question is then whether one can characterize the Sobolev spaces in the spirit of the condition (1.1). Several results have been given in this direction, for instance, a sort of converse to Theorem 1.1 due to Bagby and Ziemer [2], as well as some characterizations due to Swanson [12, 13] that are based on the Calderón-Zygmund classes (see, for example, [16, Chapter 3, p.132]). In fact, the ideas in [12] were virtually simultaneously introduced by Bourgain, Brezis, and Mironescu [3], who had shown that for $f \in L^p(\mathbb{R}^N)$ the simple condition

$$\limsup_{\epsilon \to 0} \int_{\mathbb{R}^N} \fint_{B(0,\epsilon)} \frac{|f(x+h) - f(x)|^p}{|h|^p} \, dh dx < +\infty$$

characterizes $W^{1,p}(\mathbb{R}^N)$ for $1 < p < +\infty$ (which actually holds in the context of more general mollifiers, see Section 1.2). The main idea in both [3, 12] is that instead of considering the infinitesimal behavior of averaged difference quotients one should consider the integrated infinitesimal behavior - that is, one should utilize two integrals instead of one. Notice that the case $p = 1$ is excluded, since here one obtains a characterization not of $W^{1,1}(\mathbb{R}^N)$ but $BV(\mathbb{R}^N)$, the space of functions of bounded variation. To unify the approach for all $1 \leq p < +\infty$, let us introduce the following definition.

**Definition 1.2** *A function $f : \mathbb{R}^N \to \mathbb{R}$ is said to have a first order $L^p$-Taylor approximation if $f \in L^p(\mathbb{R}^N)$ and there exists a function $v \in L^p(\mathbb{R}^N; \mathbb{R}^N)$ such that*

$$\lim_{\epsilon \to 0} \int_{\mathbb{R}^N} \fint_{B(0,\epsilon)} \frac{|f(x+h) - f(x) - v(x) \cdot h|^p}{|h|^p} \, dh dx = 0. \qquad (1.2)$$

The first result of this paper is the following theorem which asserts that functions in the Sobolev space $W^{1,p}(\mathbb{R}^N)$ possess a first order $L^p$-Taylor approximation.

**Theorem 1.3** *Let $1 \leq p < \infty$ and $f \in W^{1,p}(\mathbb{R}^N)$. Then $f$ has a first order $L^p$-Taylor approximation, and moreover, one has the stronger estimate*

$$\lim_{\epsilon \to 0} \int_{\mathbb{R}^N} \left( \fint_{B(0,\epsilon)} \frac{|f(x+h) - f(x) - \nabla f(x) \cdot h|^{pq}}{|h|^{pq}} \, dh \right)^{\frac{1}{q}} dx = 0, \qquad (1.3)$$

*where $1 \leq q \leq \frac{N}{N-p}$ if $1 \leq p < N$, $1 \leq q < \infty$ if $p = N$, and $1 \leq q \leq \infty$ if $p > N$.*

What is quite surprising is that this property - the existence of a Taylor approximation in this $L^p$ sense - in fact characterizes Sobolev functions. Precisely, we have the following theorem characterizing the Sobolev space $W^{1,p}(\mathbb{R}^N)$ in terms of the $L^p(\mathbb{R}^N)$ convergence (1.2).

**Theorem 1.4** *Let $1 \leq p < +\infty$. Then $f \in W^{1,p}(\mathbb{R}^N)$ if and only if $f$ has a first order $L^p$-Taylor approximation.*

The existence of an $L^1$-Taylor approximation is of a similar character to the assumption of Swanson in [13], since it allows one to gain some equi-integrability in an approximating sequence. If one does not make such a strong assumption, say that the limit (1.3) is finite but not zero, one obtains the following characterization of the space of functions of bounded variation.

**Theorem 1.5** *Suppose $f \in L^1(\mathbb{R}^N)$. Then $f \in BV(\mathbb{R}^N)$ if and only if there exists a function $v \in L^1(\mathbb{R}^N; \mathbb{R}^N)$ such that*

$$\limsup_{\epsilon \to 0} \int_{\mathbb{R}^N} \left( \fint_{B(0,\epsilon)} \frac{|f(x+h) - f(x) - v(x) \cdot h|^q}{|h|^q} \, dh \right)^{\frac{1}{q}} dx < +\infty \qquad (1.4)$$

*for any $1 \leq q \leq \frac{N}{N-1}$.*

In general, the limit (1.4) is bounded by a constant times $|D^s f|(\mathbb{R}^N)$, the total variation of the singular portion of the measure $Df$. However, when $q = 1$ one can explicitly compute the limit, which is given by

**Theorem 1.6** *Suppose $f \in BV(\mathbb{R}^N)$. Then writing $Df = v\mathcal{L}^N + D^s f$ where $v \in L^1(\mathbb{R}^N; \mathbb{R}^N)$ and $D^s f \perp \mathcal{L}^N$, one has*

$$\lim_{\epsilon \to 0} \int_{\mathbb{R}^N} \fint_{B(0,\epsilon)} \frac{|f(x+h) - f(x) - v(x) \cdot h|}{|h|} \, dhdx = |D^s f|(\mathbb{R}^N).$$

These results have been announced in [11] (with the exception of Theorem 1.6) with proofs of several implications, and it is our intention here to provide the complete proofs of all of the claims. However, another major point we wish to address is the fact that both Theorem 1.1 and Theorem 1.3 hold in the context of a larger framework - the use of more general approximations of the identity as seen in the work of Bourgain, Brezis, and Mironescu [3]. We now develop the necessary background for these results and several later in the paper.

## 1.2 Connections to the work of Bourgain, Brezis, and Mironescu

Let us suppose that $\Omega \subset \mathbb{R}^N$ is open, bounded, and smooth (or all of $\mathbb{R}^N$), and that $\{\rho_\epsilon\} \subset L^1(\mathbb{R}^N)$ are radial mollifiers that satisfy

$$\rho_\epsilon \geq 0, \quad \int_{\mathbb{R}^N} \rho_\epsilon(x) \, dx = 1, \qquad (1.5)$$

$$\lim_{\epsilon \to 0} \int_{|x| > \delta} \rho_\epsilon(x) \, dx = 0 \quad \text{for all } \delta > 0. \qquad (1.6)$$

Taking $\rho_\epsilon(x) = \frac{1}{|B(0,\epsilon)|} \chi_{B(0,\epsilon)}(x)$, we revert to the assumptions of the previous section. However, one has the possibility of new interesting examples, such as that of the "Gagliardo semi-norms":

$$\rho_\epsilon(x) = c_\epsilon \frac{\chi_{B(0,R)}}{|x|^{N+(\epsilon-1)r}},$$

for $c_\epsilon \to 0$ as $\epsilon \to 0$, or more generally (cf. Brezis [4])

$$\rho_\epsilon(x) = \begin{cases} 0 & \text{if } 0 < |x| \leq \epsilon, \\ a_\epsilon \psi(|x|) & \text{if } \epsilon < |x| \leq 1, \\ 0 & \text{if } 1 < |x|, \end{cases}$$

where $\psi \in L^1_{loc}(0,1)$ satisfies

$$\int_0^1 \psi(r) r^{N-1} \, dr = +\infty$$

and

$$a_\epsilon := \left( \int_\epsilon^1 \psi(r) r^{N-1} \, dr \right)^{-1} \to 0$$

as $\epsilon \to 0$.

Then the work[1] of Bourgain, Brezis, and Mironescu [3] showed two results of particular relevance here. Firstly, if $1 \leq p < +\infty$ and $f \in W^{1,p}(\Omega)$ then one has

$$\lim_{\epsilon \to 0} \int_\Omega \frac{|f(x+h) - f(x)|^p}{|h|^p} \rho_\epsilon(x-y) dy = K_{p,N} |\nabla f(x)|^p \qquad (1.7)$$

for

$$K_{p,N} := \fint_{S^{N-1}} |e_1 \cdot \sigma|^p \, d\mathcal{H}^{N-1}(\sigma).$$

Secondly, they showed that if $1 < p < +\infty$ and $f \in L^p(\Omega)$, one has $f \in W^{1,p}(\Omega)$ if and only if

$$\limsup_{\epsilon \to 0} \int_\Omega \int_\Omega \frac{|f(x+h) - f(x)|^p}{|h|^p} \rho_\epsilon(x-y) dy dx < +\infty,$$

and in that case

$$\lim_{\epsilon \to 0} \int_\Omega \int_\Omega \frac{|f(x+h) - f(x)|^p}{|h|^p} \rho_\epsilon(x-y) dy dx = K_{p,N} \int_\Omega |\nabla f(x)| \, dx. \qquad (1.8)$$

In fact, when $\rho_\epsilon(x) = \frac{1}{|B(0,\epsilon)|} \chi_{B(0,\epsilon)}(x)$ one can deduce the convergence as a consequence of $L^p$-differentiability (in an improved form allowing for $\frac{1}{|h|^p}$ in place of $\frac{1}{\epsilon^p}$, cf. Ambrosio, Fusco, and Pallara [1]), since the estimate

$$|A^p - B^p| \leq |A - B| C (1 + A^{p-1} + B^{p-1}), \qquad (1.9)$$

---

[1] We follow the convention of the subsequent paper of Ponce [14] in taking mollifiers indexed by $\epsilon$ instead of $n$.

and Hölder's inequality imply

$$\lim_{\epsilon \to 0} \left( \fint_{B(x,\epsilon)} \frac{|f(y) - f(x)|^p}{|x - y|^p} \, dy - K_{p,N}|\nabla f(x)|^p \right)$$

$$= \lim_{\epsilon \to 0} \left( \fint_{B(x,\epsilon)} \frac{|f(y) - f(x)|^p}{|x - y|^p} \, dy - \fint_{B(x,\epsilon)} \left| \nabla f(x) \cdot \frac{y - x}{|y - x|} \right|^p \, dy \right)$$

$$\leq \lim_{\epsilon \to 0} \left( \fint_{B(x,\epsilon)} \frac{|f(y) - f(x) - \nabla f(x) \cdot (y - x)|^p}{|x - y|^p} \, dy \right)^{\frac{1}{p}}$$

$$\times \left( \fint_{B(x,\epsilon)} C(1 + \frac{|f(y) - f(x)|^p}{|x - y|^p} + |\nabla f(x)|^p \, dy \right)^{\frac{1}{p'}},$$

while the difference quotient appearing in the second quantity can easily be bounded by observing that

$$\left( \fint_{B(x,\epsilon)} \frac{|f(y) - f(x)|^p}{|x - y|^p} \, dy \right)^{\frac{1}{p}}$$

$$\leq \left( \fint_{B(x,\epsilon)} \frac{|f(y) - f(x) - \nabla f(x) \cdot (y - x)|^p}{|x - y|^p} \, dy \right)^{\frac{1}{p}} + \left( \fint_{B(x,\epsilon)} \left| \nabla f(x) \cdot \frac{y - x}{|y - x|} \right| \, dy \right)^{\frac{1}{p}}.$$

The same computation follows for general mollifiers, and in this spirit, one might ask whether the stronger pointwise convergence result of Theorem 1.1 holds in this setting. Indeed, we have the following theorem to this effect.

**Theorem 1.7** Let $1 \leq p < \infty$ and $f \in W^{1,p}_{loc}(\mathbb{R}^N)$. Suppose $\rho_\epsilon$ are radial non-increasing and satisfy (1.5) and (1.6). Then for any $\eta > 0$ one has

$$\lim_{\epsilon \to 0} \int_{B(0,\eta)} \frac{|f(x + h) - f(x) - \nabla f(x) \cdot h|^{pq}}{|h|^{pq}} \rho_\epsilon(h) \, dh = 0$$

for $\mathcal{L}^N$ almost every $x \in \mathbb{R}^N$, where $1 \leq q \leq \frac{N}{N-p}$ if $1 \leq p < N$ and $1 \leq q < \infty$ if $p \geq N$.

This observation prompts us to return to the energy convergence (1.8) and perform a similar estimate. In this setting, one observes that using Minkowski's inequality for integrals implies

$$\lim_{\epsilon \to 0} \int_\Omega \left( \left( \int_\Omega \frac{|f(y) - f(x)|^{pq}}{|x - y|^{pq}} \rho_\epsilon(x - y) dy \right)^{\frac{1}{q}} - K_{pq,N}^{\frac{1}{q}} |\nabla f(x)|^p \right) dx$$

$$\leq \lim_{\epsilon \to 0} \int_\Omega \left( \left( \int_\Omega \left( \frac{|f(y) - f(x)|^p}{|x - y|^p} - \left| \nabla f(x) \cdot \frac{y - x}{|x - y|} \right|^p \right)^q \rho_\epsilon(x - y) \, dy \right)^{\frac{1}{q}} \right) dx,$$

which using again the inequality (1.9) and Hölder's inequality yields

$$\lim_{\epsilon \to 0} \int_\Omega \left( \left( \int_\Omega \frac{|f(y) - f(x)|^{pq}}{|x-y|^{pq}} \rho_\epsilon(x-y) dy \right)^{\frac{1}{q}} - K_{pq,N}^{\frac{1}{q}} |\nabla f(x)|^p \right) dx$$

$$\leq \lim_{\epsilon \to 0} \left( \int_\Omega \left( \int_\Omega \frac{|f(y) - f(x) - \nabla f(x) \cdot (y-x)|^{pq}}{|x-y|^{pq}} \rho_\epsilon(x-y) dy \right)^{\frac{1}{q}} dx \right)^{\frac{1}{p}}$$

$$\times \left( \int_\Omega \left( \int_\Omega C(1 + \frac{|f(y) - f(x)|^{pq}}{|x-y|^{pq}} + |\nabla f(x)|^{pq} \rho_\epsilon(x-y) \, dy \right)^{\frac{1}{q}} dx \right)^{\frac{1}{p'}}.$$

The same estimate as before shows that the quantities in question are finite, provided one can establish a generalization of the estimate (1.3) to other families of mollifiers. The following theorem shows that for $p > 1$ and $pq < p^*$, the Sobolev critical exponent (see Section 2), one has such a generalization.

**Theorem 1.8** Let $1 < p < \infty$ and $f \in W^{1,p}_{loc}(\mathbb{R}^N)$. Suppose $\rho_\epsilon$ are radial non-increasing and satisfy (1.5) and (1.6). Then for any $\eta > 0$ one has

$$\lim_{\epsilon \to 0} \int_K \left( \int_{B(0,\eta)} \frac{|f(x+h) - f(x) - \nabla f(x) \cdot h|^{pq}}{|h|^{pq}} \rho_\epsilon(h) \, dh \right)^{\frac{1}{q}} = 0 \qquad (1.10)$$

for all $K \subset \mathbb{R}^N$ compact, and where $1 \leq q < \frac{N}{N-p}$ if $1 < p < N$.

The degeneracy when $p = 1$ is more than an issue with $W^{1,1}$ versus $BV$, the space of functions of bounded variation. Whereas previously in the paper [10], the author and G. Leoni had given an example of a function $f \in BV(0,1)$ showing one could not expect such a convergence for $f \in BV$, we here present a counterexample to the $W^{1,1}$ case, we here in Section 5 an example (which is an adaptation a counterexample of Ponce presented by Nguyen in [15]) of a family of radial non-increasing mollifiers $\rho_\epsilon$ and a function $f \in W^{1,1}(0,1)$ for which

$$\lim_{\epsilon \to 0} \int_0^1 \left( \int_0^1 \frac{|f(x) - f(y)|^q}{|x-y|^q} \rho_\epsilon(x-y) \, dy \right)^{\frac{1}{q}} dx = +\infty$$

for any $q > 1$. In particular, this shows the restriction $p > 1$ in Theorem 1.8 is necessary. It would be interesting to understand whether for $1 < p < N$ the assumption that $pq < p^*$ is also necessary or an artifact of the proof.

The plan of the paper is as follows. In Section 2 we review our notation and recall some lemmata regarding weakly differentiable functions. In Section 3 we restrict our attention to mollification by the characteristic function, completing the proofs of Theorems 1.3, 1.4, and 1.5 announced in [11]. We also prove Theorem 1.6. In Section 4 we treat the case of general radial non-increasing mollifiers and prove the main results from Section 1.2. Finally, in Section 5 we treat the case of non-smooth domains in relation to the works [9, 10]. Here we substantiate a claim in the paper [9], provided one is willing restrict oneself to radial non-increasing mollifiers (or mollifiers that admit a family of radial non-increasing majorants - see Section 5).

## 2 Preliminaries

In what follows for $1 \leq p < N$ we denote by

$$p^* = \frac{Np}{N-p}$$

the Sobolev critical exponent and formally take $p^* = +\infty$ if $p > N$.

We use the notation

$$\fint_{B(x,r)} f(y)\, dy := \frac{1}{|B(x,r)|} \int_{B(x,r)} f(y)\, dy$$

to denote the average value of the function $f$ over the open ball centered at $x$ and of radius $r$, which we sometimes will shorten to

$$\fint f,$$

when the set of integration is implied by the context.

We utilize $S^{N-1}$ for the unit sphere in $\mathbb{R}^N$ and $\mathcal{H}^{N-1}$ for the $(N-1)$-dimensional Hausdorff measure.

For $f \in W^{1,p}_{loc}(\mathbb{R}^N)$ we write $\nabla f$ for the distributional gradient of $f$, while for $f \in BV_{loc}(\mathbb{R}^N)$ we write $Df$. Let us recall that $Df$ admits a decomposition $Df = \nabla f \mathcal{L}^N + D^s f$, where we utilize with an abuse of notation $\nabla f \in L^1_{loc}(\mathbb{R}^N; \mathbb{R}^N)$ to denote the Radon-Nikodym derivative of $Df$ with respect to the Lebesgue measure and $D^s f$ its singular measure ($D^s f \perp \mathcal{L}^N$).

We utilize the notation

$$\mathcal{M}(g)(x) := \sup_{t>0} \fint_{B(x,t)} |g(y)|\, dy,$$

to denote the Hardy-Littlewood maximal function.

We denote by $C$ a constant which may change from line to line.

Let us here record some inequalities that will be useful in the sequel. The first lemma is Poincaré's inequality for Sobolev functions on balls, a proof of which can be found as Theorem 2 on page 141 of [8].

**Lemma 2.1** *Let $1 \leq p < N$. Then there exists a constant $C = C(N, p) > 0$ such that*

$$\left( \fint_{B(x,r)} |f(y) - \fint f|^{p^*}\, dy \right)^{\frac{1}{p^*}} \leq Cr \left( \fint_{B(x,r)} |\nabla f(y)|^p\, dy \right)^{\frac{1}{p}},$$

*for all $B(x,r) \subset \mathbb{R}^N$ and $f \in W^{1,p}(B(x,r))$.*

Essentially the same proof gives the following analogous result for functions of bounded variation (which is Theorem 1 on page 189 of [8]).

**Lemma 2.2** *There exists a constant $C = C(N) > 0$ such that*

$$\left( \fint_{B(x,r)} |f(y) - \fint f|^{1^*}\, dy \right)^{\frac{1}{1^*}} \leq Cr \frac{|Df|(B(x,r))}{r^N}.$$

*for all $B(x,r) \subset \mathbb{R}^N$ and $f \in BV_{loc}(\mathbb{R}^N)$.*

We now prove the following lemma

**Lemma 2.3** *Suppose $f \in W^{1,p}_{loc}(\mathbb{R}^N)$ for some $1 \leq p < \infty$, and that $1 \leq q < +\infty$ is such that $pq \leq p^*$. Then there exists a $C = C(p,q,N) > 0$ such that for all $t > 0$*

$$\frac{1}{t^{N+pq}} \int_{B(0,t)} |f(x+h) - f(x) - \nabla f(x) \cdot h|^{pq}\, dh$$

$$\leq C \left( \fint_{B(0,t)} |\nabla f(x+h) - \nabla f(x)|^s\, dh \right)^{\frac{pq}{s}}$$

$$+ C \left( \frac{1}{t} \fint_{B(0,t)} |f(x+z) - f(x) - \nabla f(x) \cdot z|\, dz \right)^{pq},$$

*where $1 \leq s \leq p$ is such that $s^* = pq$.*

**Proof.** We begin with the estimate

$$\frac{1}{t^{N+pq}} \int_{B(0,t)} |f(x+h) - f(x) - \nabla f(x) \cdot h|^{pq}\, dh$$

$$\leq \frac{C}{t^{N+pq}} \int_{B(0,t)} |f(x+h) - \fint f - \nabla f(x) \cdot h|^{pq}\, dh$$

$$+ \frac{C}{t^{N+pq}} \int_{B(0,t)} |\fint_{B(0,t)} f(x+z) - f(x) - \nabla f(x) \cdot z\, dz|^{pq}\, dh$$

where we have added and subtracted $\fint f(x+z)\, dz$ (where the integral is taken over $B(0,t)$) and used the fact that $\int \nabla f(x) \cdot z\, dz = 0$. The second term agrees with what we aim to show in this lemma, so it remains to show that

$$\frac{1}{t^{N+pq}} \int_{B(0,t)} |f(x+h) - \fint f - \nabla f(x) \cdot h|^{pq}\, dh \leq C' \left( \fint_{B(0,t)} |\nabla f(x+h) - \nabla f(x)|^s\, dh \right)^{\frac{pq}{s}}$$

and the result is demonstrated. However, utilizing the fact that $1 \leq s \leq p$, we have that $f \in W^{1,s}_{loc}(\mathbb{R}^N)$, and so applying Lemma 2.1 to

$$v(h) := f(x+h) - \fint f - \nabla f(x) \cdot h$$

(note that $f \in W^{1,s}_{loc}(\mathbb{R}^N)$ implies $v \in W^{1,s}(B(0,R))$ for every $0 < R < \infty$) we obtain

$$\int_{B(0,t)} |f(x+h) - \fint f - \nabla f(x) \cdot h|^{pq}\, dh$$

$$\leq C \left( \int_{B(0,t)} |\nabla f(x+h) - \nabla f(x)|^s\, dh \right)^{\frac{pq}{s}},$$

since $\fint v = 0$ and $\nabla v(h) = \nabla f(x+h) - \nabla f(x)$. Then dividing by $t^{N+pq}$, we obtain the desired inequality. ∎

Analogously we have

**Lemma 2.4** *Suppose $f \in BV_{loc}(\mathbb{R}^N)$ that $1 \leq q \leq 1^*$. Then there exists a $C = C(q, N) > 0$ such that for all $t > 0$*

$$\frac{1}{t^{N+q}} \int_{B(0,t)} |f(x+h) - f(x) - \nabla f(x) \cdot h|^q \, dh$$

$$\leq C \left( \fint_{B(0,t)} |\nabla f(x+h) - \nabla f(x)| \, dh + |D^s f|(B(0,t)) \right)^q$$

$$+ C \left( \frac{1}{t} \fint_{B(0,t)} |f(x+z) - f(x) - \nabla f(x) \cdot z| \, dz \right)^q.$$

**Proof.** As before

$$\frac{1}{t^{N+q}} \int_{B(0,t)} |f(x+h) - f(x) - \nabla f(x \cdot) h|^q \, dh$$

$$\leq \frac{C}{t^{N+q}} \int_{B(0,t)} |f(x+h) - \fint f - \nabla f(x) \cdot h|^q \, dh$$

$$+ \frac{C}{t^{N+q}} \int_{B(0,t)} |\fint_{B(0,t)} f(x+z) - f(x) - \nabla f(x) \cdot z \, dz|^q \, dh,$$

though here instead we have the estimate

$$\frac{1}{t^{N+q}} \int_{B(0,t)} |f(x+h) - \fint f - \nabla f(x) \cdot h|^q \, dh = \frac{C}{t^q} \fint_{B(0,t)} |f(x+h) - \fint f - \nabla f(x) \cdot h|^q \, dh$$

$$\leq \frac{C}{t^q} \left( \fint_{B(0,t)} |f(x+h) - \fint f - \nabla f(x) \cdot h|^{1^*} \, dh \right)^{\frac{q}{1^*}}$$

$$= \left( \frac{C}{t^{1^*}} \fint_{B(0,t)} |f(x+h) - \fint f - \nabla f(x) \cdot h|^{1^*} \, dh \right)^{\frac{q}{1^*}}.$$

Now, applying Lemma 2.2 to $v(h) := f(x+h) - \fint f - \nabla f(x) \cdot h$, we have

$$\left( \frac{C}{t^{1^*}} \fint_{B(0,t)} |f(x+h) - \fint f - \nabla f(x) \cdot h|^{1^*} \, dh \right)^{\frac{1}{1^*}} \leq C \fint_{B(0,t)} |\nabla f(x+h) - \nabla f(x)| dx + |D^s f|(B(0,t)),$$

and this implies the desired result. ∎

Finally, when $p > N$ we have Morrey's inequality for Sobolev functions with a precise estimate for the exponent of Hölder continuity as given in the following Lemma (which is stated and proved as Theorem 3 on page 189 of [8]).

**Lemma 2.5** *Let $N < p < +\infty$. There exists a $C = C(p, N) > 0$ such that*

$$|f(z) - f(y)| \leq Cr \left( \fint_{B(x,r)} |\nabla f(w)|^p \, dw \right)^{\frac{1}{p}}$$

*for all $B(x, r) \subset \mathbb{R}^N$, $f \in W^{1,p}(B(x,r))$, and $\mathcal{L}^N$ almost every $y, z \in B(x,r)$.*

# 3  $L^p$-Taylor approximations and $W^{1,p}$

Let us recall that in [11] we gave a proof of Theorem 1.3 in the regime $1 \leq p < N$, Theorem 1.4 and the necessity of $f \in BV(\mathbb{R}^N)$ supposing that the lim sup is bounded in Theorem 1.5. Therefore we here complete the arguments to Theorem 1.3 in the regime $p \geq N$ and Theorem 1.5. We also prove Theorem 1.6.

We begin with the proof of Theorem 1.3 in the regime $p \geq N$.

**Proof.** Let $p \geq N$ and $f \in W^{1,p}(\mathbb{R}^N)$. If $q < +\infty$, then the result follows by exactly the same argument in [11] observing that although Lemma 2.3 was only stated for $p < N$ in [11], we here proved it for all $1 \leq p < +\infty$ and any $1 \leq q < +\infty$. Thus it remains to treat the case $q = +\infty$. The estimate for Morrey's inequality given in Lemma 2.5 implies that for any $N < s < p$ one has

$$\frac{|f(x+h) - f(x) - \nabla f(x) \cdot h|^s}{|h|^s} \leq C \fint_{B(x,|h|)} |\nabla f(z) - \nabla f(x)|^s \, dz$$

and so

$$\sup_{h \in B(0,\epsilon)} \frac{|f(x+h) - f(x) - \nabla f(x) \cdot h|^p}{|h|^p} \leq C \sup_{h \in B(0,\epsilon)} \left( \fint_{B(x,|h|)} |\nabla f(z) - \nabla f(x)|^s \, dz \right)^{\frac{p}{s}}.$$

Now, the right hand side tends to zero pointwise as $\epsilon \to 0$, while

$$\sup_{h \in B(0,\epsilon)} \left( \fint_{B(x,|h|)} |\nabla f(z) - \nabla f(x)|^s \, dz \right)^{\frac{p}{s}} \leq \mathcal{M}(|\nabla f|^s)^{\frac{p}{s}}(x) + |\nabla f(x)|^p,$$

and so the result follows from Lebesgue's dominated convergence theorem and the boundedness of $\mathcal{M} : L^r(\mathbb{R}^N) \to L^r(\mathbb{R}^N)$ for all $1 < p \leq +\infty$ ∎

We now give a proof of Theorem 1.5, arguing as in [11], mutatis mutandis.

**Proof.** As argued in [11], the finiteness of the limit implies $f \in BV(\mathbb{R}^N)$, so it suffices to show that for $f \in BV(\mathbb{R}^N)$ we can find a bound for the limit.

For any $0 < \epsilon < 1$, we expand the integrand on concentric rings

$$\fint_{B(0,\epsilon)} \frac{|f(x+h) - f(x) - \nabla f(x) \cdot h|^q}{|h|^q} \, dh$$

$$= \sum_{i=0}^{\infty} \frac{1}{\epsilon^N |B(0,1)|} \int_{B(0,\frac{\epsilon}{2^i}) \setminus B(0,\frac{\epsilon}{2^{i+1}})} \frac{|f(x+h) - f(x) - \nabla f(x) \cdot h|^q}{|h|^q} \, dh.$$

We now make estimates for $i \in \mathbb{N}$ fixed. We have

$$\frac{1}{\epsilon^N} \int_{B(0,\frac{\epsilon}{2^i}) \setminus B(0,\frac{\epsilon}{2^{i+1}})} \frac{|f(x+h) - f(x) - \nabla f(x) \cdot h|^q}{|h|^q} \, dh$$

$$\leq \frac{1}{\epsilon^N} \left(\frac{\epsilon}{2^{i+1}}\right)^{-q} \int_{B(0,\frac{\epsilon}{2^i}) \setminus B(0,\frac{\epsilon}{2^{i+1}})} |f(x+h) - f(x) - \nabla f(x) \cdot h|^q \, dh$$

$$\leq \frac{2^q}{2^{iN}} \left(\frac{\epsilon}{2^i}\right)^{-N-q} \int_{B(0,\frac{\epsilon}{2^i})} |f(x+h) - f(x) - \nabla f(x) \cdot h|^q \, dh.$$

Here, in place of Lemma 2.3 we utilize Lemma 2.4 to deduce that

$$\left(\frac{\epsilon}{2^i}\right)^{-N-q} \int_{B(0,\frac{\epsilon}{2^i})} |f(x+h) - f(x) - \nabla f(x) \cdot h|^q \, dh$$

$$\leq C \left( \fint_{B(0,\frac{\epsilon}{2^i})} |\nabla f(x+h) - \nabla f(x)|^p \, dh + 2^{iN} |D^s f|(\overline{B(x, 1/2^i)}) \right)^q$$

$$+ C \left( \int_0^1 \fint_{B(0,t\frac{\epsilon}{2^i})} |\nabla f(x+tz) - \nabla f(x)| \, dz + \left(\frac{2^i}{t}\right)^N |D^s f|(\overline{B(x, t/2^i)}) dt \right)^q.$$

Therefore, summing in $i$ and applying the basic inequality $(\sum_i |a_i|)^{\frac{1}{q}} \leq \sum_i |a_i|^{\frac{1}{q}}$ (which follows from subadditivity of the function $s \mapsto s^{\frac{1}{q}}$), we have

$$\left( \fint_{B(0,\epsilon)} \frac{|f(x+h) - f(x) - \nabla f(x) \cdot h|^q}{|h|^q} \, dh \right)^{\frac{1}{q}}$$

$$\leq C \sum_{i=0}^\infty \left(\frac{1}{2^i}\right)^{N/q} \fint_{B(0,\frac{\epsilon}{2^i})} |\nabla f(x+h) - \nabla f(x)| \, dh + 2^{iN} |D^s f|(\overline{B(x, 1/2^i)})$$

$$+ C \sum_{i=0}^\infty \left(\frac{1}{2^i}\right)^{N/q} \int_0^1 \fint_{B(0,t\frac{\epsilon}{2^i})} |\nabla f(x+tz) - \nabla f(x)| \, dz + \left(\frac{2^i}{t}\right)^N |D^s f|(\overline{B(x, t/2^i)}) dt.$$

Integrating the preceding inequality over $x \in \mathbb{R}^N$ and making use of Tonelli's theorem we obtain

$$\int_{\mathbb{R}^N} \left( \fint_{B(0,\epsilon)} \frac{|f(x+h) - f(x) - \nabla f(x) \cdot h|^q}{|h|^q} \, dh \right)^{\frac{1}{q}} dx$$

$$\leq C \sum_{i=0}^\infty \left(\frac{1}{2^i}\right)^{N/q} \fint_{B(0,\frac{\epsilon}{2^i})} \int_{\mathbb{R}^N} |\nabla f(x+h) - \nabla f(x)| \, dx dh + \int_{\mathbb{R}^N} 2^{iN} |D^s f|(\overline{B(x, 1/2^i)}) dx$$

$$+ C \sum_{i=0}^\infty \left(\frac{1}{2^i}\right)^{N/q} \int_0^1 \left( \fint_{B(0,t\frac{\epsilon}{2^i})} \int_{\mathbb{R}^N} |\nabla f(x+tz) - \nabla f(x)| \, dz + \int_{\mathbb{R}^N} \left(\frac{2^i}{t}\right)^N |D^s f|(\overline{B(x, t/2^i)}) \, dx \right) dt.$$

However, if $h, z \in B(0, \epsilon)$ and $t \in (0,1)$ we have

$$\max \left\{ \int_{\mathbb{R}^N} |\nabla f(x+h) - \nabla f(x)|^p \, dx, \int_{\mathbb{R}^N} |\nabla f(x+tz) - \nabla f(x)|^p \, dx \right\}$$

$$\leq \sup_{w \in B(0,\epsilon)} \int_{\mathbb{R}^N} |\nabla f(x+w) - \nabla f(x)|^p \, dx,$$

and observe that this bound is independent of $i \in \mathbb{N}$. Thus,

$$\int_{\mathbb{R}^N} \left( \fint_{B(0,\epsilon)} \frac{|f(x+h) - f(x) - \nabla f(x) h|^q}{|h|^q} \, dh \right)^{\frac{1}{q}} dx$$

$$\leq C \sum_{i=0}^\infty \left(\frac{1}{2^i}\right)^{N/q} \sup_{w \in B(0,\epsilon)} \int_{\mathbb{R}^N} |\nabla f(x+w) - \nabla f(x)|^p \, dx$$

$$+ C \sum_{i=0}^\infty \left(\frac{1}{2^i}\right)^{N/q} \left( \int_{\mathbb{R}^N} 2^{iN} |D^s f|(\overline{B(x, 1/2^i)}) dx + \int_0^1 \int_{\mathbb{R}^N} \left(\frac{2^i}{t}\right)^N |D^s f|(\overline{B(x, t/2^i)}) \, dx dt \right).$$

Finally we observe that

$$\int_0^1 \int_{\mathbb{R}^N} \left(\frac{2^i}{t}\right)^N |D^s f|(\overline{B(x,t/2^i)}) \, dxdt = \int_0^1 \left(\frac{2^i}{t}\right)^N \int_{B(0,t\frac{1}{2^i})} \int_{\mathbb{R}^N} d|D^s f|(x-y) \, dydt$$
$$\leq C|D^s f|(\mathbb{R}^N)$$

and similarly for

$$\int_{\mathbb{R}^N} 2^{iN} |D^s f|(\overline{B(x,1/2^i)}) dx.$$

Thus, as the infinite series is summable, sending $\epsilon \to 0$ and using continuity of translation in $L^1(\mathbb{R}^N)$ we obtain

$$\lim_{\epsilon \to 0} \int_{\mathbb{R}^N} \left( \fint_{B(0,\epsilon)} \frac{|f(x+h) - f(x) - \nabla f(x)h|^q}{|h|^q} \, dh \right)^{\frac{1}{q}} dx \leq C|D^s f|(\mathbb{R}^N).$$

∎

We now calculate the explicit limit of the $L^1$-Taylor approximation for a function $f \in BV(\mathbb{R}^N)$, proving the claim of Theorem 1.6. Observe here that the fact that $q=1$ means we do not need to use the Sobolev embedding theorem in the form of Poincaré's inequality and therefore the computation is cleaner.

**Proof.** We first argue the upper bound

$$\limsup_{\epsilon \to 0} \int_{\mathbb{R}^N} \fint_{B(0,\epsilon)} \frac{|f(x+h) - f(x) - \nabla f(x) \cdot h|}{|h|} \, dhdx \leq K_{1,N} |D^s f|(\mathbb{R}^N) \tag{3.1}$$

Fix $\psi \in C_c^\infty(B(0,1))$ radial, non-negative and $\int \psi = 1$ and define $\psi_\delta(x) := \frac{1}{\delta^N} \psi(\frac{x}{\delta})$. For $f \in BV(\mathbb{R}^N)$ we denote by $f_\delta$ the convolution $f * \psi_\delta$. Then if $Df$ is the measure derivative of $f$, one has $Df_\delta = (\nabla f)_\delta + (D^s f)_\delta$, and so the triangle inequality and Tonelli's theorem imply

$$\int_{\mathbb{R}^N} \fint_{B(0,\epsilon)} \frac{|f_\delta(x+h) - f_\delta(x) - (\nabla f)_\delta(x) \cdot h|}{|h|} \, dhdx$$
$$\leq \int_0^1 \int_{\mathbb{R}^N} \fint_{B(0,\epsilon)} |(\nabla f_\delta(x+sh) - (\nabla f)_\delta(x)) \cdot \frac{h}{|h|}| \, dhdxds$$
$$\leq \int_0^1 \int_{\mathbb{R}^N} \fint_{B(0,\epsilon)} |(\nabla f)_\delta(x+sh) - (\nabla f)_\delta(x)| \, dhdxds$$
$$+ \int_0^1 \int_{\mathbb{R}^N} \fint_{B(0,\epsilon)} \left|(D^s f)_\delta(x+sh) \cdot \frac{h}{|h|}\right| \, dhdxds$$
$$= \int_0^1 \fint_{B(0,\epsilon)} \int_{\mathbb{R}^N} |(\nabla f)_\delta(x+sh) - \nabla f_\delta(x)| \, dhdxds + \fint_{B(0,\epsilon)} \int_{\mathbb{R}^N} \left|(D^s f)_\delta(y) \cdot \frac{h}{|h|}\right| \, dydh.$$

Now

$$\fint_{B(0,\epsilon)} \int_{\mathbb{R}^N} \left| (D^s f)_\delta(y) \cdot \frac{h}{|h|} \right| \, dydh = \frac{1}{\epsilon^N |B(0,1)|} \int_0^\epsilon \int_{\partial B(0,t)} \int_{\mathbb{R}^N} |(D^s f)_\delta(y) \cdot \sigma| \, dyd\mathcal{H}^{N-1}(\sigma)dt$$

$$= \frac{1}{\epsilon^N |B(0,1)|} \int_0^\epsilon \int_{\partial B(0,t)} \int_{\mathbb{R}^N} |(D^s f)_\delta(y) \cdot \sigma| \, dyd\mathcal{H}^{N-1}(\sigma)dt$$

$$= \frac{1}{\epsilon^N |B(0,1)|} \int_0^\epsilon t^{N-1} \, dt \int_{S^{N-1}} \int_{\mathbb{R}^N} |(D^s f)_\delta(y) \cdot \sigma| \, dyd\mathcal{H}^{N-1}(\sigma)$$

$$= K_{1,N} |(D^s f)_\delta|(\mathbb{R}^N)$$

for

$$K_{1,N} = \fint_{S^{N-1}} |e_1 \cdot \sigma| \, d\mathcal{H}^{N-1}(\sigma)$$

as in Section 1.2. Therefore sending $\delta \to 0^+$, Fatou's lemma, Lebesgue's dominated convergence theorem, the strict convergence of the measures $(D^s f)_\delta \to D^s f$ yields

$$\int_{\mathbb{R}^N} \fint_{B(0,\epsilon)} \frac{|f(x+h) - f(x) - \nabla f(x) \cdot h|}{|h|} \, dhdx$$

$$\leq \int_0^1 \fint_{B(0,\epsilon)} \int_{\mathbb{R}^N} |\nabla f(x+sh) - \nabla f(x)| \, dxdhds + K_{1,N} |D^s f|(\mathbb{R}^N)$$

Now taking the limsup as $\epsilon \to 0$ the first term vanishes as argued in the previous theorems (by continuity of translation on $L^1(\mathbb{R}^N)$), and we thus obtain (3.1).

To prove the lower bound, we simply estimate

$$\int_{\mathbb{R}^N} \fint_{B(0,\epsilon)} \frac{|f(x+h - f(x)|}{|h|} \, dhdx - K_{1,N} \int_{\mathbb{R}^N} |\nabla f(y)| \, dy$$

$$= \int_{\mathbb{R}^N} \fint_{B(0,\epsilon)} \frac{|f(x+h - f(x)|}{|h|} \, dhdx - \fint_{B(0,\epsilon)} \int_{\mathbb{R}^N} \left| \nabla f(y) \cdot \frac{h}{|h|} \right| \, dydh$$

$$\leq \int_{\mathbb{R}^N} \fint_{B(0,\epsilon)} \frac{|f(x+h - f(x) - \nabla f(x) \cdot h|}{|h|} \, dhdx,$$

and taking the liminf as $\epsilon \to 0$, an application of the result of Juan Dávila [7] that

$$\lim_{\epsilon \to 0} \int_{\mathbb{R}^N} \fint_{B(0,\epsilon)} \frac{|f(x+h - f(x)|}{|h|} \, dhdx = K_{1,N} |Df|(\mathbb{R}^N),$$

and the decomposition

$$|Df|(\mathbb{R}^N) = |D^s f|(\mathbb{R}^N) + \int_{\mathbb{R}^N} |\nabla f(y)| \, dy$$

yields the desired lower bound. ∎

## 4 General Mollifiers

We will here give proofs of the results asserted in Section 1.2. We begin with Theorem 1.7.

**Proof of Theorem 1.7.** We will show that

$$\lim_{\epsilon \to 0} \left( \int_{B(0,\eta)} \frac{|f(x+h) - f(x) - \nabla f(x)h|^{pq}}{|h|^{pq}} \rho_\epsilon(h) \, dh \right)^{\frac{1}{pq}} = 0$$

for $p$-Lebesgue points of $\nabla f$, and hence $\mathcal{L}^N$ almost every $x \in \mathbb{R}^N$. Let us first recall that as a monotone decreasing function of $|h|$, $\rho_\epsilon$ has at most countably many jump discontinuities, and possesses inner and outer limits for every value of $|h|$, and therefore its suffices to prove the result in the case where $\rho_\epsilon$ is inner continuous (which always dominates $\rho_\epsilon$).

In order to apply Lemma 2.3, we expand the desired quantity on concentric rings to obtain

$$\int_{B(0,\eta)} \frac{|f(x+h) - f(x) - \nabla f(x)h|^{pq}}{|h|^{pq}} \rho_\epsilon(h) \, dh$$

$$= \sum_{i=0}^\infty \int_{B(0,\frac{\eta}{2^i}) \setminus B(0,\frac{\eta}{2^{i+1}})} \frac{|f(x+h) - f(x) - \nabla f(x)h|^{pq}}{|h|^{pq}} \rho_\epsilon(h) \, dh$$

$$\leq \sum_{i=0}^\infty \rho_\epsilon\left(\frac{\eta}{2^{i+1}}\right) \left(\frac{\eta}{2^{i+1}}\right)^{-pq} \int_{B(0,\frac{\eta}{2^i})} |f(x+h) - f(x) - \nabla f(x)h|^{pq} \, dh$$

$$= 2^{pq+2N} \sum_{i=0}^\infty \rho_\epsilon\left(\frac{\eta}{2^{i+1}}\right) \left(\frac{\eta}{2^{i+2}}\right)^N \left(\frac{\eta}{2^i}\right)^{-pq-N} \int_{B(0,\frac{\eta}{2^i})} |f(x+h) - f(x) - \nabla f(x)h|^{pq} \, dh.$$

Then applying the estimate from Lemma 2.3, we have

$$\int_{B(0,\eta)} \frac{|f(x+h) - f(x) - \nabla f(x)h|^{pq}}{|h|^{pq}} \rho_\epsilon(h) \, dh \qquad (4.1)$$

$$\leq C \sum_{i=0}^\infty \rho_\epsilon\left(\frac{\eta}{2^{i+1}}\right) \left(\frac{\eta}{2^{i+2}}\right)^N \max\left\{F_1^i(x), F_2^i(x)\right\}, \qquad (4.2)$$

where

$$F_1^i(x) := \left( \fint_{B(0,\frac{\eta}{2^i})} |\nabla f(x+h) - \nabla f(x)|^s \, dh \right)^{\frac{pq}{s}}$$

$$F_2^i(x) := \left( \left(\frac{\eta}{2^i}\right)^{-N-1} \int_{B(0,\frac{\eta}{2^i})} |f(z) - f(x) - \nabla f(x)z| \, dz \right)^{pq}.$$

Now, the Lebesgue differentiation theorem implies that for almost every $x \in \mathbb{R}^N$, $F_1^i, F_2^i$ tend to zero for large $i$, while they are bounded for $i$ small (and the bound may depend on $x$ at this point). Thus, we have that

$$\int_{B(0,\eta)} \frac{|f(x+h) - f(x) - \nabla f(x)h|^{pq}}{|h|^{pq}} \rho_\epsilon(h) \, dh$$

$$\leq \delta(x) \sum_{i=k(\delta)}^\infty \rho_\epsilon\left(\frac{\eta}{2^{i+1}}\right) \left(\frac{\eta}{2^{i+2}}\right)^N$$

$$+ C \sum_{i=0}^{k(\delta)} \rho_\epsilon\left(\frac{\eta}{2^{i+1}}\right) \left(\frac{\eta}{2^{i+2}}\right)^N.$$

We claim that the first sum is bounded by a constant which is independent of epsilon and delta, while the second tends to zero as $\epsilon \to 0$. If this is the case, then first sending $\epsilon \to 0$ and then sending $\delta \to 0$, the result is demonstrated. However, observe that since $\rho_\epsilon$ is decreasing and $t^{N-1}$ is increasing,

$$\rho_\epsilon\left(\frac{\eta}{2^{i+1}}\right)\left(\frac{\eta}{2^{i+2}}\right)^{N-1} \leq \inf_{t \in \left[\frac{\eta}{2^{i+2}}, \frac{\eta}{2^{i+1}}\right]} \rho_\epsilon(t) t^{N-1},$$

and therefore

$$\rho_\epsilon\left(\frac{\eta}{2^{i+1}}\right)\left(\frac{\eta}{2^{i+2}}\right)^{N} \leq \int_{\frac{\eta}{2^{i+2}}}^{\frac{\eta}{2^{i+1}}} \rho_\epsilon(t) t^{N-1}\, dt.$$

This then implies that

$$\sum_{i=k(\delta)}^{\infty} \rho_\epsilon\left(\frac{\eta}{2^{i+1}}\right)\left(\frac{\eta}{2^{i+2}}\right)^{N} \leq \int_0^{\frac{\eta}{k(\delta)+1}} \rho_\epsilon(t)\, t^{N-1}\, dt$$

$$\leq 1,$$

and the first sum is bounded independent of epsilon and delta, as claimed. For the demonstration that the second tends to zero as $\epsilon \to 0$, we utilize the fact that $\rho_\epsilon$ are decreasing and

$$\lim_{\epsilon \to 0} \int_\delta^\infty \rho_\epsilon(t) t^{N-1} dt = 0$$

for all $\delta > 0$. This implies that $\rho_\epsilon(t_0) \to 0$ pointwise for almost every $t_0 > 0$ and so by monotonicity of $\rho_\epsilon$, it suffices to bound the second sum by

$$\sum_{i=0}^{k(\delta)} \rho_\epsilon\left(\frac{\eta}{2^{i+1}}\right)\left(\frac{\eta}{2^{i+2}}\right)^{N} \leq \rho_\epsilon(t_0) \eta^N$$

for some $t_0 < \frac{\eta}{2^{k(\delta)+1}}$ such that $\rho_\epsilon(t_0) \to 0$. ∎

We now prove Theorem 1.8 concerning the improved first order $L^p$-Taylor approximation estimate for general mollifiers.

**Proof of Theorem 1.8.** As we demonstrated the pointwise convergence

$$\left(\int_{B(0,\eta)} \frac{|f(x+h) - f(x) - \nabla f(x)h|^{pq}}{|h|^{pq}} \rho_\epsilon(h)\, dh\right)^{\frac{1}{pq}} \to 0$$

as $\epsilon \to 0$ for $\mathcal{L}^N$ almost every $x \in \mathbb{R}^N$, it remains to show that there is a function $g \in L^1_{loc}(\mathbb{R}^N)$ such that

$$\left(\int_{B(0,\eta)} \frac{|f(x+h) - f(x) - \nabla f(x)h|^{pq}}{|h|^{pq}} \rho_\epsilon(h)\, dh\right)^{\frac{1}{q}} \leq g(x)$$

and the result is demonstrated. Recall that we showed

$$\left(\int_{B(0,\eta)} \frac{|f(x+h) - f(x) - \nabla f(x)h|^{pq}}{|h|^{pq}} \rho_\epsilon(h)\, dh\right)^{\frac{1}{q}}$$

$$\leq C\left(\sum_{i=0}^{\infty} \rho_\epsilon\left(\frac{\eta}{2^{i+1}}\right)\left(\frac{\eta}{2^{i+2}}\right)^{N} \max\left\{F_1^i(x), F_2^i(x)\right\}\right)^{\frac{1}{q}}$$

where

$$F_1^i(x) = \left(\fint_{B(0,\frac{\eta}{2^i})} |\nabla f(x+h) - \nabla f(x)|^s \, dh\right)^{\frac{pq}{s}}$$

$$F_2^i(x) = \left(\left(\frac{\eta}{2^i}\right)^{-N-1} \int_{B(0,\frac{\eta}{2^i})} |f(z) - f(x) - \nabla f(x)z| \, dz\right)^{pq}.$$

Now, we can estimate uniformly in $i$

$$F_1^i(x) \leq \mathcal{M}_\eta(|\nabla f|^s(x))^{\frac{pq}{s}} + |\nabla f(x)|^{pq}$$
$$F_2^i(x) \leq \mathcal{M}_\eta(|\nabla f|(x))^{pq} + |\nabla f(x)|^{pq},$$

where $\mathcal{M}_\eta$ is the (restricted) Hardy-Littlewood maximal function,

$$\mathcal{M}_\eta(h)(x) := \sup_{0<r<\eta} \fint_{B(x,r)} |h(z)| \, dz$$

so that

$$\left(\int_{B(0,\eta)} \frac{|f(x+h) - f(x) - \nabla f(x)h|^{pq}}{|h|^{pq}} \rho_\epsilon(h) \, dh\right)^{\frac{1}{q}} \leq g(x),$$

where we have defined

$$g(x) := C\left(\mathcal{M}(|\nabla f|^s(x))^{\frac{p}{s}} + |\nabla f(x)|^p + \mathcal{M}(|\nabla f|(x))^p\right) \sup_{\epsilon>0} \left(\sum_{i=0}^\infty \rho_\epsilon\left(\frac{\eta}{2^{i+1}}\right)\left(\frac{\eta}{2^{i+2}}\right)^N\right)^{\frac{1}{q}}.$$

It therefore remains to show that $g \in L^1_{loc}(\mathbb{R}^N)$. Now, in analogy to proof of the previous theorem, one can show that

$$\sup_{\epsilon>0} \sum_{i=0}^\infty \rho_\epsilon\left(\frac{\eta}{2^{i+1}}\right)\left(\frac{\eta}{2^{i+2}}\right)^N \leq \sup_{\epsilon>0} \int_0^\infty \rho_\epsilon(t) t^{N-1} \, dt$$
$$= 1,$$

and so it remains to show that $\mathcal{M}(|\nabla f|(x))^p, \mathcal{M}(|\nabla f|^s(x))^{\frac{p}{s}}, \in L^1_{loc}(\mathbb{R}^N)$. For both of these functions, we simply use the fact that the restricted Hardy-Littlewood maximal function is a bounded map from $L^r_{loc}(\mathbb{R}^N)$ to $L^r_{loc}(\mathbb{R}^N)$ for $r > 1$. In the first case, we use the fact that $p > 1$ and $\nabla f \in L^p_{loc}(\mathbb{R}^N)$, while for the second, we utilize fact that $\frac{p}{s} > 1$ and $|\nabla f|^s \in L^{\frac{p}{s}}_{loc}(\mathbb{R}^N)$. The fact that $\frac{p}{s} > 1$ is a simple computation when $p \geq N$, while follows from the assumption $q < \frac{N}{N-p}$ when $1 < p < N$. ∎

## 5 Application to Characterizations of Sobolev spaces

In the recent paper [9], we made claim to a new necessary and sufficient condition for functions in the Sobolev space $W^{1,p}(\Omega)$ where $\Omega \subset \mathbb{R}^N$ is an arbitrary open set. The sufficiency of the condition was established, though a problem with the argument of the necessity required that a different argument be used to establish this direction. One of the main reasons we developed these arguments and these tools in this paper is for this purpose, to complete the claimed result in regard to Sobolev functions. Let us first recall some of the notation and required prerequisites developed in [9].

Given $\{\rho_\epsilon\}$ possibly non-radial, we assume that there exist $\{v_i\}_{i=1}^N \subset \mathbb{R}^N$, linearly independent, and a $\delta > 0$ such that

$$\liminf_{\epsilon \to 0} \int_{C_\delta(v_i)} \rho_\epsilon(x)\, dx > 0 \text{ for all } i = 1, \ldots, N, \tag{5.1}$$

where

$$C_\delta(v) := \left\{ w \in \mathbb{R}^N \backslash \{0\} : \frac{v}{|v|} \cdot \frac{w}{|w|} > 1 - \delta \right\}.$$

Then if one defines the measures, for $E \subset S^{N-1}$,

$$\mu_\epsilon(E) := \int_E \int_0^\infty \rho_\epsilon(t\sigma) t^{N-1}\, dt d\mathcal{H}^{N-1}(\sigma),$$

it is shown in [9] that up to a subsequence, $\mu_{\epsilon_j} \overset{*}{\rightharpoonup} \mu$ weakly-star in $\left(C_0(S^{N-1})\right)'$ and one has that there exists $\alpha > 0$ such that

$$\int_{S^{N-1}} |v \cdot \sigma|\, d\mu(\sigma) \geq \alpha |v|$$

for all $v \in \mathbb{R}^N$.

Also recall the definition of $\Omega_\lambda$, defined by

$$\Omega_\lambda := \{ x \in \Omega : |x| < \frac{1}{\lambda}, \operatorname{dist}(x, \partial\Omega) > \lambda \}.$$

In a recent paper, the author and G. Leoni had stated the following theorem, whose proof is correct in the case $q = 1$.

**Theorem 5.1** *Let $\Omega \subset \mathbb{R}^N$ be open, let $\{\rho_\epsilon\} \subset L^1(\mathbb{R}^N)$ satisfy (1.5), (1.6), and (5.1), let $1 \leq p < \infty$ and $1 \leq q < \infty$, with $1 \leq q \leq \frac{N}{N-p}$ if $p < N$, and let $f \in L^1_{\operatorname{loc}}(\Omega)$. Then $f \in W^{1,1}_{\operatorname{loc}}(\Omega)$ and $\nabla f \in L^p(\Omega; \mathbb{R}^N)$ if and only if*

$$\lim_{\lambda \to 0} \limsup_{\epsilon \to 0} \int_{\Omega_\lambda} \left( \int_{\Omega_\lambda} \left( \frac{|f(x) - f(y)|^p}{|x-y|^p} \right)^q \rho_\epsilon(x - y)\, dy \right)^{\frac{1}{q}} dx < +\infty. \tag{5.2}$$

*Moreover, there exist a subsequence $\{\epsilon_j\}$ and a probability measure $\mu \in M(S^{N-1})$ such that*

$$\lim_{\lambda \to 0} \lim_{j \to \infty} \int_{\Omega_\lambda} \left( \int_{\Omega_\lambda} \left( \frac{|f(x) - f(y)|^p}{|x-y|^p} \right)^q \rho_{\epsilon_j}(x - y)\, dy \right)^{\frac{1}{q}} dx$$

$$= \int_\Omega \left( \int_{S^{N-1}} (|\nabla f(x) \cdot \sigma|^p)^q\, d\mu(\sigma) \right)^{\frac{1}{q}} dx.$$

We here aim to show that if $f \in W^{1,p}(\Omega)$, then

$$\lim_{\lambda \to 0} \limsup_{\epsilon \to 0} \int_{\Omega_\lambda} \left( \int_{\Omega_\lambda} \left( \frac{|f(x) - f(y)|^p}{|x - y|^p} \right)^q \rho_\epsilon(x - y) \, dy \right)^{\frac{1}{q}} dx < +\infty,$$

as well as that there exist a subsequence $\{\epsilon_j\}$ and a probability measure $\mu \in M(S^{N-1})$ such that

$$\lim_{\lambda \to 0} \lim_{j \to \infty} \int_{\Omega_\lambda} \left( \int_{\Omega_\lambda} \left( \frac{|f(x) - f(y)|^p}{|x - y|^p} \right)^q \rho_{\epsilon_j}(x - y) \, dy \right)^{\frac{1}{q}} dx$$
$$= \int_\Omega \left( \int_{S^{N-1}} (|\nabla f(x) \cdot \sigma|^p)^q \, d\mu(\sigma) \right)^{\frac{1}{q}} dx.$$

We do so under some minor additional hypothesis. First, we suppose that $1 < p$ and $pq < p^*$. Secondly we suppose the family $\rho_\epsilon$ admits a family of non-increasing radial majorant $\hat{\rho}_\epsilon$ with similar properties. (In particular, if $\rho_\epsilon$ are radial non-increasing this is satisfied, and so this hypothesis is not empty), where non-increasing radial majorants are defined as follows.

**Definition 5.2** *We say that a family of mollifiers $\{\rho_\epsilon\}$ admits a family of non-increasing radial majorants $\{\hat{\rho}_\epsilon\}$ if there exists $\eta > 0$ such that*

$$\rho_\epsilon(h) \leq \hat{\rho}_\epsilon(h)$$

*for all $|h| \leq \eta$ and $\hat{\rho}_\epsilon(t)$ are non-increasing and satisfy*

$$\int_0^\infty \hat{\rho}_\epsilon(t) t^{N-1} \, dt = C$$
$$\lim_{\epsilon \to 0} \int_\delta^\infty \hat{\rho}_\epsilon(t) t^{N-1} \, dt = 0$$

*for every $\delta > 0$.*

Assuming Theorem 1.8, we will demonstrate the correct upper bound needed for Theorem 1.5 from our paper, recalling that we only need to consider the truncated version

$$\int_{\Omega_\lambda} \left( \int_{B(x,\eta)} \left( \frac{|f(x) - f(y)|^p}{|x - y|^p} \right)^q \rho_\epsilon(x - y) \, dy \right)^{\frac{1}{q}} dx,$$

for $\eta > 0$ which is arbitrarily small (and therefore small enough to apply the radial majorant hypothesis), since the limit as $\epsilon \to 0$ of

$$\int_{\Omega_\lambda} \left( \int_{\Omega_\lambda \setminus B(x,\eta)} \left( \frac{|f(x) - f(y)|^p}{|x - y|^p} \right)^q \rho_\epsilon(x - y) \, dy \right)^{\frac{1}{q}} dx$$

was shown to vanish.

**Proof.** If we define

$$A_\epsilon(x) := \left( \int_{B(0,\eta)} \frac{|f(x+h) - f(x)|^{pq}}{|h|^{pq}} \rho_\epsilon(h) \, dh \right)^{\frac{1}{pq}}$$

$$B_\epsilon(x) := \left( \int_{B(0,\eta)} \left| \nabla f(x) \cdot \frac{h}{|h|} \right|^{pq} \rho_\epsilon(h) \, dh \right)^{\frac{1}{pq}},$$

then our aim is to show that

$$\lim_{\epsilon \to 0} \int_{\Omega_\lambda} |A_\epsilon^p(x) - B_\epsilon^p(x)| \, dx = 0.$$

If this is the case, we will have completed the proof of the theorem, since the upper bound

$$\int_{\Omega_\lambda} B_\epsilon^p(x) \, dx \leq \int_{\Omega_\lambda} |\nabla f(x)|^p \, dx$$
$$\leq \int_{\Omega} |\nabla f(x)|^p \, dx,$$

implies that the desired quantity is bounded as we take the limit in $\epsilon$ and $\lambda$, while along the subsequence $\epsilon_j$ we have

$$\lim_{j \to \infty} \int_{\Omega_\lambda} B_{\epsilon_j}^p(x) \, dx = \int_{\Omega_\lambda} \left( \int_{S^{N-1}} |\nabla f(x) \cdot \sigma|^{pq} \, d\mu(\sigma) \right)^{\frac{1}{q}} dx.$$

However, by convexity of the function $t \mapsto t^p$ and the bounds we have on $A_\epsilon$ and $B_\epsilon$, it suffices to show that

$$\lim_{\epsilon \to 0} \int_{\Omega_\lambda} |A_\epsilon(x) - B_\epsilon(x)|^p \, dx = 0, \tag{5.3}$$

since we may estimate

$$\int_{\Omega_\lambda} |A_\epsilon^p(x) - B_\epsilon^p(x)| \, dx$$
$$\leq C \left( \int_{\Omega_\lambda} |A_\epsilon(x) - B_\epsilon(x)|^p \, dx \right)^{\frac{1}{p}} \left( \int_{\Omega_\lambda} 1 + A_\epsilon^p(x) + B_\epsilon^p(x) \, dx \right)^{\frac{1}{p'}},$$

and where $p'$ is conjugate to $p$ in the sense that $\frac{1}{p} + \frac{1}{p'} = 1$ (and notice that when $p = 1$ we do not even need this estimate). Notice also that we have shown the necessary bounds for $B_\epsilon^p$, while $A_\epsilon^p \leq C \left( |A_\epsilon - B_\epsilon|^p + B_\epsilon^p \right)$. Thus, we will establish (5.3). We have that

$$|A_\epsilon(x) - B_\epsilon(x)| \leq \left( \int_{B(0,\eta)} \frac{|f(x+h) - f(x) - \nabla f(x)h|^{pq}}{|h|^{pq}} \rho_\epsilon(h) \, dh \right)^{\frac{1}{pq}}$$
$$\leq \left( \int_{B(0,\eta)} \frac{|f(x+h) - f(x) - \nabla f(x)h|^{pq}}{|h|^{pq}} \hat{\rho}_\epsilon(h) \, dh \right)^{\frac{1}{pq}},$$

where we have utilized the assumption that the family $\rho_\epsilon$ admits a family of non-increasing radial majorants. However, now an application of Theorem 1.8 implies that the quantity on the right-hand-side tends to zero in $L^p_{loc}$, and hence we conclude the convergence (5.3) and the result is demonstrated. ∎

Now the preceding proof precluded the possibility $p = 1$, and so one might wonder whether another proof might be possible in that regime, supposing one has a genuine $W^{1,1}$ function. The following counterexample shows that this is not the case, so that the assumption $p > 1$ is sharp.

**Counterexample 5.3** *We adapt a counterexample due to Ponce which can be found in the paper of Nguyen [15]. Let $q > 1$, define $\delta_n := \frac{1}{2^n}$ and choose $\epsilon_n \to 0$ such that*

$$\delta_{n+1}\epsilon_n^{\frac{1}{q}-1}(\delta_{n+1})^{\epsilon_n/q} \geq n.$$

*Then if we define $m_n := \frac{\delta_n + \delta_{n+1}}{2}$, we claim*

$$\delta_{n+1} \int_{\delta_{n+1}}^{m_n} \left( \int_{m_n}^{\delta_n} \frac{c\epsilon_n}{|x-y|^{q+1-\epsilon_n}} \, dy \right)^{\frac{1}{q}} dx \geq \tilde{c}n.$$

*To see this, note that evaluating the inner integral one has*

$$\int_{m_n}^{\delta_n} \frac{1}{|x-y|^{q+1-\epsilon_n}} \, dy = (q - \epsilon_n) \left( \frac{1}{|x-m_n|^{q-\epsilon_n}} - \frac{1}{|x-\delta_n|^{q-\epsilon_n}} \right),$$

*while the inequality $(a-b)^{\frac{1}{q}} \geq a^{\frac{1}{q}} - b^{\frac{1}{q}}$ yields*

$$\delta_{n+1} \int_{\delta_{n+1}}^{m_n} \left( \int_{m_n}^{\delta_n} \frac{c\epsilon_n}{|x-y|^{q+1-\epsilon_n}} \, dy \right)^{\frac{1}{q}} dx$$

$$\geq \delta_{n+1}(c\epsilon_n)^{\frac{1}{q}}(q-\epsilon_n) \int_{\delta_{n+1}}^{m_n} \frac{1}{|x-m_n|^{1-\epsilon_n/q}} - \frac{1}{|x-\delta_n|^{1-\epsilon_n/q}} \, dx$$

$$= \delta_{n+1}(c\epsilon_n)^{\frac{1}{q}}(q-\epsilon_n)\frac{q}{\epsilon_n} \left( (m_n - \delta_{n+1})^{\epsilon_n/q} - (\delta_n - m_n)^{\epsilon_n/q} + (\delta_n - \delta_{n+1})^{\epsilon_n/q} \right)$$

$$= \delta_{n+1}(c\epsilon_n)^{\frac{1}{q}}(q-\epsilon_n)\frac{q}{\epsilon_n}(\delta_{n+1})^{\epsilon_n/q}$$

$$\geq \tilde{c}\delta_{n+1}\epsilon_n^{\frac{1}{q}-1}(\delta_{n+1})^{\epsilon_n/q}.$$

Thus, we may find $\eta_n > 0$ such that

$$\delta_{n+1} \int_{\delta_{n+1}}^{m_n-\eta_n} \left( \int_{m_n+\eta_n}^{\delta_n} \frac{c\epsilon}{|x-y|^{q+1-\epsilon}} \, dy \right)^{\frac{1}{q}} dx \geq \frac{1}{2}\tilde{c}n. \tag{5.4}$$

Now, we define

$$f(x) := \begin{cases} \delta_{n+1} : x \in [\delta_{n+1}, m_n - \eta_n] \\ \delta_n : x \in [m_n + \eta_n, \delta_n] \end{cases}$$

and affine on $[m_n - \eta_n, m_n + \eta_n]$. Then $f \in W^{1,1}(0,1)$, since $f$ is piecewise affine and increasing, while

$$\lim_{n \to \infty} \int_0^1 \left( \int_0^1 \frac{|f(x)-f(y)|^q}{|x-y|^q} \frac{c\epsilon_n}{|x-y|^{1-\epsilon_n}} \, dy \right)^{\frac{1}{q}} dx = +\infty.$$

To see this, note that for each $n \in \mathbb{N}$, we have

$$\int_0^1 \left( \int_0^1 \frac{|f(x)-f(y)|^q}{|x-y|^q} \frac{c\epsilon_n}{|x-y|^{1-\epsilon_n}} \, dy \right)^{\frac{1}{q}} dx \geq \delta_{n+1} \int_{\delta_{n+1}}^{m_n-\eta_n} \left( \int_{m_n+\eta_n}^{\delta_n} \frac{c\epsilon_n}{|x-y|^{q+1-\epsilon_n}} \, dy \right)^{\frac{1}{q}} dx,$$

and the result follows from the inequality (5.4)

# Acknowledgements

The author would like to thank Augusto Ponce for his discussions regarding this work and encouragement during the project. The author is supported by the Taiwan Ministry of Science and Technology under research grant MOST 103-2115-M-009-016-MY2.